\documentclass[11pt]{article}

\usepackage{makeidx}
\usepackage{geometry}                
\usepackage{graphicx}
\usepackage{amssymb}
\usepackage{epstopdf}
\usepackage{lscape}
\usepackage{pdflscape}
\usepackage{verbatim}

\usepackage[gen]{eurosym}

\newcommand{\pf}{\vs \noindent {\it Proof:} \quad}

\newcommand{\thm}{\vspace{.25in} \noindent {\bf Theorem.} \quad }

\newcommand{\C}{\mbox{\msbm{C}}}

\newcommand{\R}{\mbox{\msbm{R}}}

\newcommand{\be}{\begin{equation}}
\newcommand{\ee}{\end{equation}}
\newcommand{\ba}{\begin{array}}
\newcommand{\ea}{\end{array}}
\newcommand{\bea}{\begin{eqnarray}}
\newcommand{\eea}{\end{eqnarray}}
\newcommand{\bean}{\begin{eqnarray*}}
\newcommand{\eean}{\end{eqnarray*}}

\newcommand{\vs}{\vspace{.25in}}

\font\msbm=msbm10

\def\picture #1 by #2 (#3){
  \vbox to #2{
    \hrule width #1 height 0pt depth 0pt
    \vfill
    \special{picture #3} 
    }
  }

\def\scaledpicture #1 by #2 (#3 scaled #4){{
  \dimen0=#1 \dimen1=#2
  \divide\dimen0 by 1000 \multiply\dimen0 by #4
  \divide\dimen1 by 1000 \multiply\dimen1 by #4
  \picture \dimen0 by \dimen1 (#3 scaled #4)}
  }

\newcommand{\bc}{\begin{center}}
\newcommand{\ec}{\end{center}}
\newcommand{\aster}{\begin{center} *** \end{center}}

\begin{document}

\title{
An Algebraic Realization of the Taylor-Socolar Aperiodic Monotilings
}
\author{H. L. Resnikoff\footnote{Resnikoff Innovations LLC; howard@resnikoff.com}}
\date{20141211}         

\maketitle

\begin{abstract}
The first aperiodic monotiling, introduced by Taylor, was based on  a  trapezoidal prototile equipped with 14 distinct decorations. A presentation of the closely related  Taylor-Socolar aperiodic monotiling  is based on a hexagonal prototile equipped with 7 decorations. This paper gives  decoration-free algebraic descriptions equivalent to each of these presentations. It also shows how the monotilings and  Taylor triangles pattern that characterizes the aperiodicity can be obtained from just one algebraic equation. \\  

{\sc Keywords:} {\it  Aperiodic tilings, decoration,  inflation,  monotile, substitution tilings, Socolar-Taylor monotile, Taylor monotile, Taylor triangles.}
\end{abstract}

\aster

A recurrent theme in the history of mathematics is the interplay between geometry and algebra. Descartes took the first steps, which made it possible to express as algebraic equations the properties of curves that had been studied by Greek mathematicians. In the opposite direction, it became possible to draw `pictures' of equations. This led to new methods for proving theorems and new insights into the relationships between algebraic and geometric structures.\footnote{The early developments are described in chapter 8 of \cite{MIC}.} This interplay continues as a major theme of contemporary mathematics. The present paper is presented in this spirit. Its purpose is to provide a simple algebraic characterization of the monotilings recently introduced into the theory of aperiodic plane tilings by Taylor, and by Taylor and Socolar. Their results were obtained by ingenious but essentially conventional geometrical methods familiar in the theory of tilings. 

In 2009 Joan M. Taylor \cite{Taylor 2010} discovered a monotile that produces only aperiodic tilings of the plane. This monotile -- a trapezoid obtained by dividing a regular hexagon into two congruent parts -- was equipped with 14 distinct `decorations' -- curves drawn on the trapezoid --  that govern which tile edges are permitted to be adjacent in a tiling. Some of the decorations are related by reflection so that a decorated tile and its mirror image are not equivalent with respect to translations and rotations.  The geometrical substitution rule requires that the curves be continuous across tile boundaries: thus, the curves  begin and end at infinity, or they are closed in compact subsets of the plane. The  curves form a characteristic pattern that can be described as the union of a straight line and an infinite collection of nested equilateral triangles. We refer to this as the {\it Taylor triangles} pattern. Triangles of increasing and unbounded size appear as larger regions of the plane are examined. This implies that the tiling is aperiodic (with respect to the decoration rules) because a periodic pattern cannot contain a system of similar figures of increasing and unlimited size. Each application of Taylor's original geometrical substitution rule generates a similar trapezoid of twice the linear dimensions; thus it is `self-similar' and after $n$ `inflations', the inflated trapezoid is tiled by $4^n$ tiles drawn from the original marked set of 14 congruent tile classes. If the center of inflation is chosen as an interior point of the trapezoid, then the  tiling covers the plane.

Socolar and Taylor followed up this discovery with a related construction that dropped the requirement of self-similarity. They showed  that a regular hexagon could be decorated in 7 different ways  to generate an aperiodic tiling by an analogous geometrical substitution \cite{Socolar+Taylor 2012, Socolar+Taylor 2011, Socolar+Taylor 2010}. The decorations lead to the same set of curves and the same argument is used to prove aperiodicity.  In addition, and more important, they constructed a geometrical shape -- a modified hexagon called  the {\it monotile} -- that encodes the rules in its shape alone and generates the same tiling. See \cite{Baake+Gaehler+Grimm (2012), Goodman-Strauss,Gruenbaum+Shephard} for additional background, and \cite{Socolar talk 2013} for potential applications to physics. 

\aster

In this paper we translate the geometrical substitutions into algebraic equations. This is both a conceptual and a calculational simplification, and brings the theory of tilings into closer contact with other parts of mathematics. 

Consider the plane equipped with the euclidean metric and the measure it induces. Two subsets of the plane are {\it essentially disjoint} if the measure of their intersection is 0.  Two subsets of the plane are {\it essentially identical} if the measures of both sets and their intersection are equal. A {\it tiling} of the plane by congruent copies of a tile $T$ is a cover by translated, rotated and possibly reflected copies of $T$ that are essentially disjoint. An {\it over-tiling} of the plane by congruent copies of a tile $T$ is a cover by translated, rotated and possibly reflected copies of $T$ that are essentially disjoint or essentially identical. An over-tiling may contain more than one copy of a tile in some positions. Notice that the allowed motions constitute all the transformations of the plane that preserve the metric.

Identify $\R^2$ with the complex plane so that points are complex numbers: $(x,y) \mapsto z= x  + i y$. Denote complex conjugation by an overline: $z \rightarrow \overline{z}$.

\begin{figure}[t]
\begin{center}
\caption{Socolar-Taylor monotile with supplementary black bar.}
\label{barred S-T monotile}
\includegraphics[width=2in]{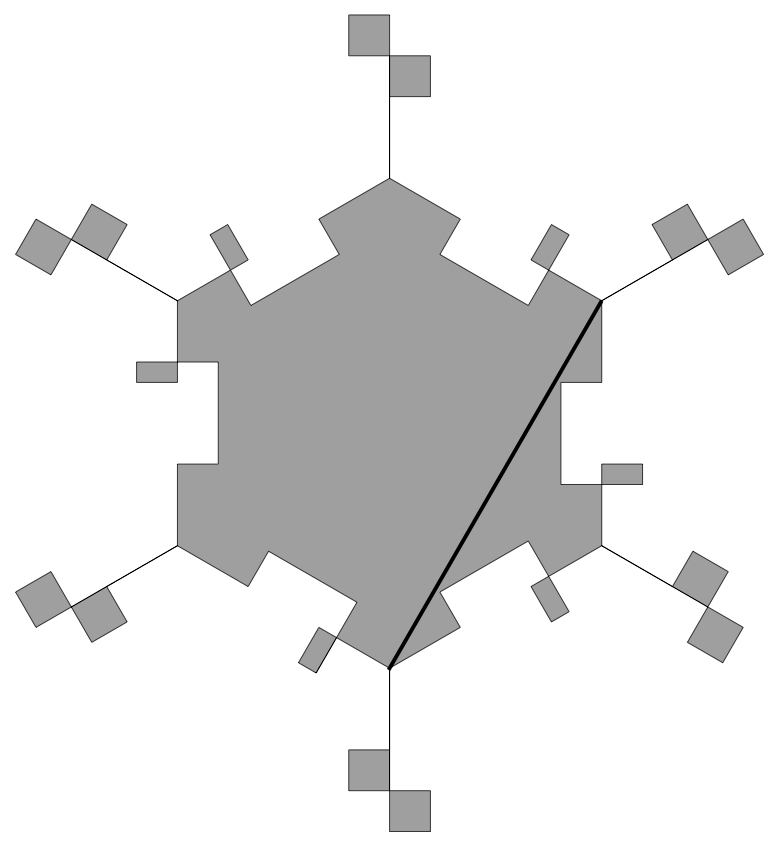}
\end{center}
\end{figure}

The Socolar-Taylor (hereafter ``S-T") monotile is a subset of the plane. The interior of the monotile consists of 19 connected components:  a large central region,  12 large  `flags', and 6 small flags.  A version of the monotile is displayed in fig.\ref{barred S-T monotile}. A  thick black chord, to be used in the later discussion,  has been drawn on it connecting two vertices of  the underlying hexagon that are separated by a third vertex. The thin black lines (of width 0) shown in the figure are only for purposes of visualization.

\begin{figure}[t]
\begin{center}
\caption{The partial aperiodic tiling $2^2 R_3 (2)$ produced by eq(\ref{HLR S-T set eqs}) showing the Taylor triangles and the underlying configuration of monotiles with black bands.}
\label{HLR 2 digit monotiling}
\includegraphics[width=3.5in]{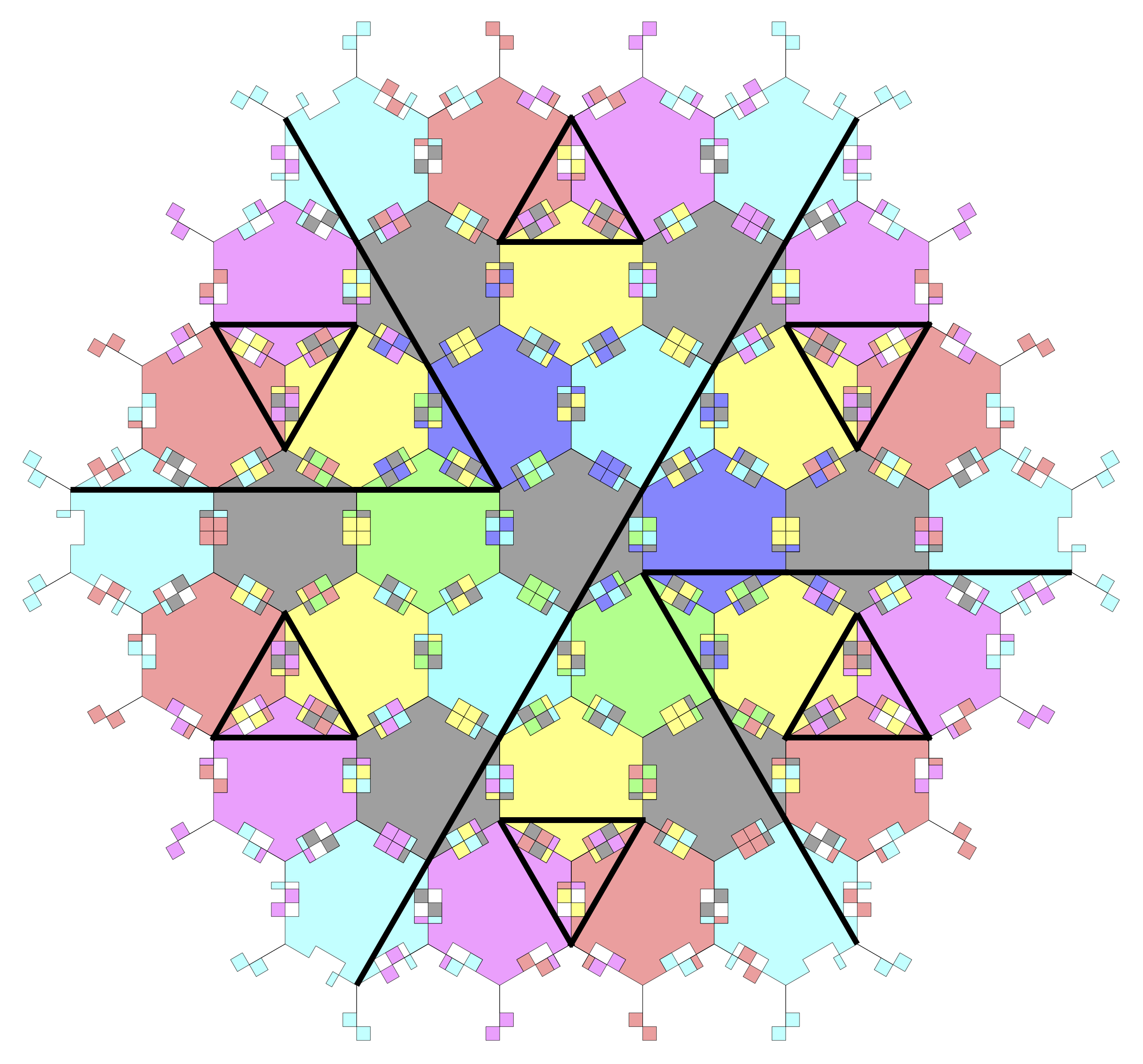}
\end{center}
\end{figure}

S-T use 7 copies of a regular hexagon decorated with distinct decorations  to construct a system of geometrical substitutions that generate the same aperiodic tiling. The substitutions have been presented in a convenient diagrammatic form by Frettl\"{o}h  in \cite{Frettloeh 2014}, where  the base figure is a regular hexagon distinguished by 7 colors.  In addition to the colors, each hexagon is decorated with the same pattern of lines and dots.   In tilings, one of the decorated tile types occurs with twice the frequency of the others. In the figures below we color this one, which corresponds to the type  `C' trapezoids in Taylor's original paper and to the third hexagonal tile in Frettl\"{o}h's  diagram \cite{Frettloeh 2014},  gray.  The other 6 tiles, taken in the same order as in Taylor and Frettl\"{o}h, are colored red, yellow, green, cyan, blue, and magenta;\footnote{These are not the colors used in \cite{Frettloeh 2014}.} cp. fig.\ref{HLR 2 digit monotiling} which shows the second iteration of the tiling process for the `C' tile, here labeled $R_3$. 

Denote the  decorated monotiles $R_k, \, 1 \leq k \leq 7$. We will consider the limit sets of recursions.   Each recursion is of the form  
$$\rho\,R_k (n) = \bigcup_{j=1}^7  \left( c_{k,j} + u_{k,j} \check{ R}_j(n-1) \right)$$
 where $\rho=2$, $\check{ R}_j(n-1)$ is either $R_j(n-1)$ or $\overline{R_j(n-1)}$ and  the  $c_{j,k}$, $u_{j,k}$ are constants. 

\thm Set  $\omega = e^{i \pi /3}$. The following system of equations describes an over-tiling of the plane that is equivalent to the  Socolar-Taylor tiling. The equations are to be interpreted as relations between limit sets of a recursion.
\be
\ba{rcl}
\rho \,  R_1 &=& 	R_3 \cup  
		\left( 1+ \omega^5 R_7 \right) \cup
		\left( \omega + R_4   \right) \cup
		\left( \omega^2 + \omega^3 \overline{R_7}   \right) \cup \\
	&& 	\qquad \left( \omega^3 + \omega^2 R_1   \right) \cup
		\left( \omega^4 + \omega^5  \overline{R_4}     \right) \cup
		\left( \omega^5 +   \overline{R_1}     \right) \\
\rho  \, R_2 &=&	  
		R_3 \cup 
		\left( 1 + \omega^4 \overline{R_2}   \right) \cup
		 \left( \omega + R_6   \right) \cup
		 \left( \omega^2 +   \omega^3 \overline{R_7}  \right)  \cup  \\
	&& 	\qquad
		 \left( \omega^3 + \omega^2 R_1 \right) \cup
		 \left( \omega^4 + \omega^5  \overline{R_6}     \right)  \cup
		\left( \omega^5 +   \omega R_7   \right) \\
\rho  \, R_3 &=&	R_3 \cup  
		\left( 1 + \omega^4  \overline{R_6}     \right) \cup	
		\left( \omega + \omega^5 \overline{R_5}   \right) \cup 
		\left( \omega^2 +    \omega^4 R_6 \right) \\
	&& 	\qquad 
		\left(  \omega^3 +\omega  \overline{ R_4}  \right) \cup
		\left( \omega^4 +   R_5   \right) \cup
		\left( \omega^5 +  \omega R_4   \right) \cup \\
\rho  \, R_4 &=&	 
		R_3 \cup 
		\left( 1 + \omega^5 R_2  \right) \cup  
		\left( \omega +  R_4 \right)  \cup 
		\left( \omega^2 + \omega^3  \overline{R_2}     \right) \cup \\
	&& 	\qquad 
		\left(  \omega^3 +\omega  \overline{ R_1}  \right) \cup
		\left( \omega^4 +    \omega^5 \overline{R_5} \right) 
		\left( \omega^5 +  \omega R_1   \right) \cup  \\
\rho  \, R_5 &=&
		R_3  \cup  
		\left(  1 + \omega^4  \overline{ R_2}  \right) \cup
		 \left( \omega + R_5   \right) \cup
		 \left( \omega^2 +  \omega^4  R_2   \right) \cup  \\
	&& 	\qquad 
		\left( \omega^3 +   \omega \overline{R_7} \right)  \cup
		\left(   \omega^4 + \omega^5  \overline{ R_5}  \right) \cup
		\left( \omega^5 +    \omega  R_7 \right) \\
\rho  \, R_6 &=&	R_3  \cup
		\left( 1 +  \omega^4  \overline{R_2}   \right) \cup
		\left( \omega +   R_6  \right)  \cup
		\left(   \omega^2 + \omega^3  \overline{ R_2}  \right) \cup  \\
	&& 	\qquad 
		\left(   \omega^3 + \omega  \overline{ R_1}  \right) \cup
		 \left( \omega^4 +   \omega^5 \overline{R_5} \right)  \cup
		\left( \omega^5 +    \omega  R_7 \right) \\
\rho  \, R_7 &=&	
		R_3  \cup		
		\left(   1 + \omega^5 R_2 \right) \cup   
		 \left( \omega +  R_4  \right)  \cup
		 \left( \omega^2 +    \omega^3  \overline{R_7}  \right)   \cup \\
	&& 	\qquad 
		\left(   \omega^3 + \omega^2   R_1  \right) \cup
		 \left( \omega^4 +    \omega^5 \overline{R_6}  \right)  \cup
		\left(  \omega^5 +  \omega R_1  \right) 
\ea
\label{HLR S-T set eqs}
\ee

\pf  The proof consists of a transcription into the algebraic language of the complex number field of the  geometric substitutions, e.g. given in \cite{Frettloeh 2014}, for the seven decorated hexagons. 

From the figure, the sets  appearing in the union appear to be essentially disjoint. They are not. Indeed, eq(\ref{HLR S-T set eqs}) implies relationships for the measure of the sets that appear in it.  Denote the measure of a measurable set $S \subset \C$ by $m(S)$. If $S_1$ and $S_2$ are essentially disjoint, then  $m (c_1 S_1 \cup c_2 S_2) = |c_1|^2 m(S_1) + |c_2|^2 m(S_2)$. Were the tiles on the right side of eq(\ref{HLR S-T set eqs})  essentially disjoint, then, since  the $R_j$ have the same measure and  $| \omega |=1$,  it would follow that $4 m(R_j ) = 7 m(R_j )$, a contradiction. Thus the equations describe an over-tiling wherein 3 of every 7 tiles overlap perfectly. $\Box$ \\

The  recursion $R_j(n-1) \mapsto R_j(n)$  implied by the theorem can be initiated by selecting  appropriate $R_j (0)$. The regular (undecorated) hexagon or the S-T monotile are instructive choices. If, for instance, the monotile is selected and decorated with an extraneous black bar as shown in fig.\ref{barred S-T monotile} merely to illustrate how the Taylor triangles originate, after two iterations one finds the  image shown above in fig.\ref{HLR 2 digit monotiling}. Note that from this algebraic perspective,  the monotile plays no role in generating either the Taylor triangles nor the general arrangement in the figure. All of the information is coded in the equations.

The theorem specifies  a recursive inflationary construction of the S-T aperiodic tiling. The tiles are either essentially disjoint or essentially identical (This simply means that we are not concerned with overlap of tile boundaries.).   As $n$ increases, the tiling covers increasing larger portion of the plane; in the limit,  $\R^2$ is tiled.  The flags fit into and fill  the gaps along each edge of the underlying hexagons.  The pattern of Taylor triangles is the same as the one generated by the geometric substitution method of Socolar and Taylor.  Socolar and Taylor prove aperiodicity by noting that pattern of the Taylor triangles  contains triangles of ever increasing size. From the point of view of the geometric description, this is the essential element in demonstrating aperiodicity. Their argument applies equally well here.

This algebraic analysis eliminates the geometric ingredients.  Everything  is contained in the algebraic description of metric-preserving transformations in the plane.

\begin{figure}[t]
\begin{center}
\caption{Marked regular hexagon prototile.}
\label{marked hex prototile}
\includegraphics[width=1.5in]{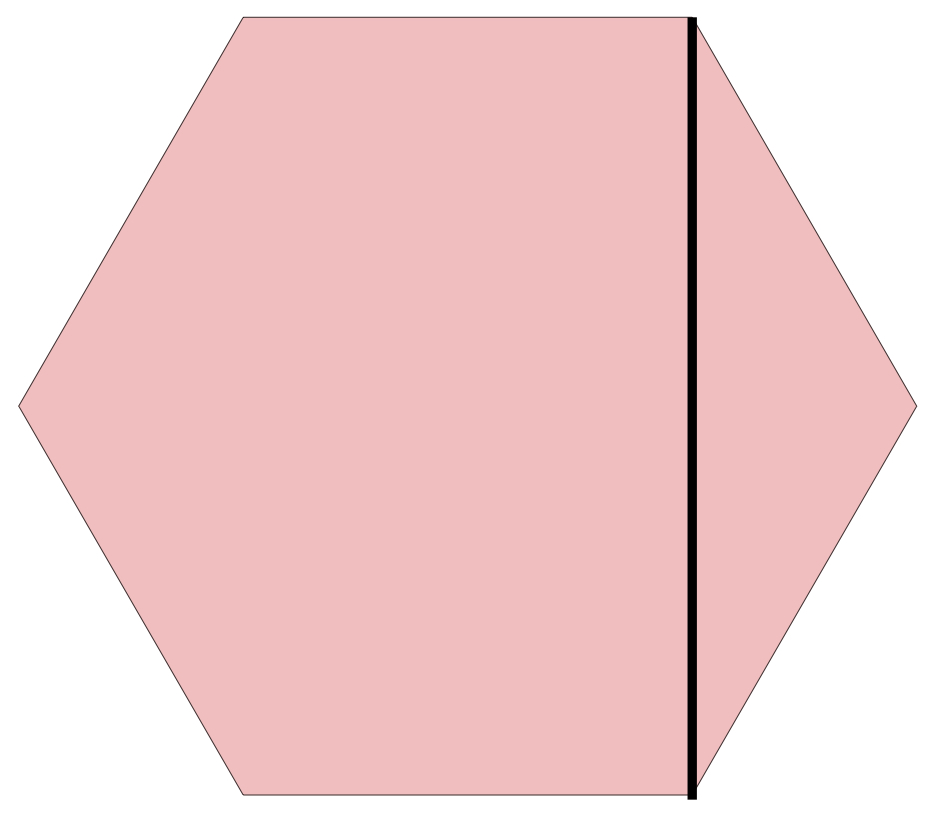}
\end{center}
\end{figure}

We can take a further step in this direction. Examination of the 7 cases in eq(\ref{HLR S-T set eqs}) reveals that the equations are all the same except for the particular arrangement of tiles $R_j$ and the  occurrence of complex conjugation. The $R_j$ differ only in color so consider one tile $R$ with one chord drawn on it, as shown in fig.\ref{marked hex prototile}. Algebraically, this amounts to substituting $R_j \rightarrow R$ and $\overline{R_j} \rightarrow \omega^5 R$. After this substitution, the seven recursion equations are identical:

\be
\ba{lcl}
\rho \,  R(n) &=& 	R(n-1) \cup  
		\left( 1+ \omega^5 R(n-1) \right) \cup
		\left( \omega + R(n-1)   \right) \cup
		\left( \omega^2 + \omega^4 R(n-1)  \right) \cup \\
	&& 	\qquad \left( \omega^3 + \omega^2 R(n-1)   \right) \cup
		\left( \omega^4 + R(n-1)     \right) \cup
		\left( \omega^5 +  \omega  R(n-1)   \right) \\
\ea
\label{remainder set R_n-eqs}
\ee

\begin{figure}[h]
\begin{center}
\caption{Taylor triangle pattern emerges from 3-iterations of the algebraic identity eq(\ref{2R limit}). The hexagon background enables one to see the how the black chords rotate.}
\label{3digit hex tiling}
\includegraphics[width=2.5in]{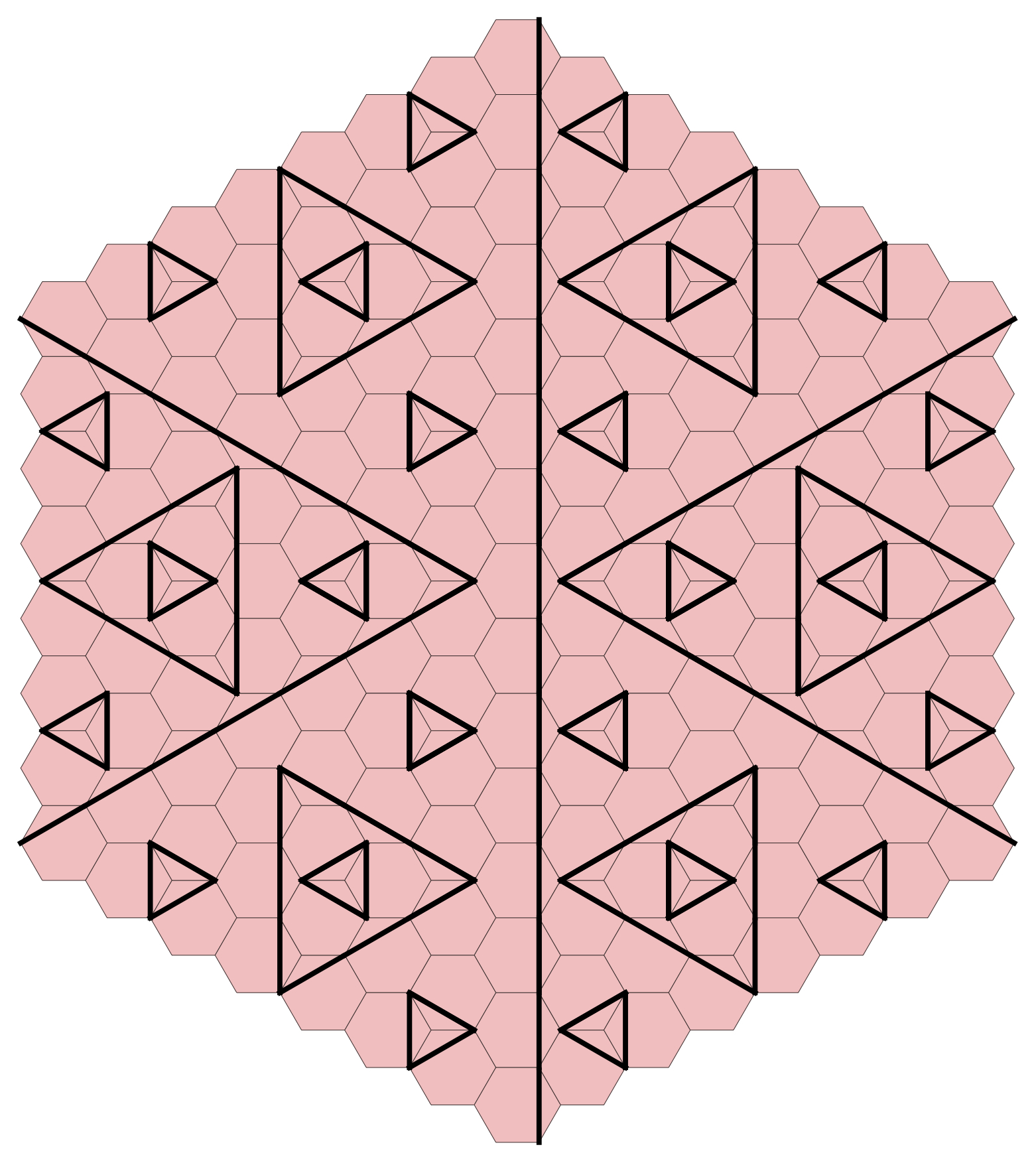}
\end{center}
\end{figure}

Figure \ref{3digit hex tiling},\footnote{The diagram has been rotated to emphasize the symmetries.} the third iteration of this unique  recursion, shows how the Taylor triangles are produced by this substitution that no longer involves reflections, i.e. complex conjugation.  The limiting version of this recursion is

\be
\ba{lcl}
\rho  \, R &=& R \cup  
		\left( 1+ \omega^5 R \right) \cup
		\left( \omega + R   \right) \cup
		\left( \omega^2 + \omega^4 R    \right) \cup   \left( \omega^3 + \omega^2 R   \right) \cup
		\left( \omega^4 + R    \right) \cup
		\left( \omega^5 +  \omega R    \right)
\\
\ea
\label{2R limit}
\ee

Again, the recursions in eq(\ref{remainder set R_n-eqs}, \ref{2R limit}) are not self-similar and they describe an over-tiling. The corresponding tiling is obtained by identifying overlaid tiles.

\thm Equation\,(\ref{2R limit}) describes an aperiodic tiling.

\pf As before, the tiling is aperiodic because the system of Taylor triangles contains triangles of arbitrarily large size. $\Box$\\

\aster

Just as the use of many different decorations may not seem entirely in the spirit of the search for a monotile, and just as a monotile defined solely by its shape may seem somewhat unsatisfactory if the shape consists of disconnected subsets of the plane\footnote{Taylor and Socolar raised this point.}, so too some readers may find the above algebraic equations somewhat unsatisfactory because they produce over-tilings rather than ordinary tilings. This deficiency can be eliminated by returning to Taylor's original trapezoids with 14 varieties of decorations. It has the additional advantage of supplying a simpler algebraic presentation of the Taylor monotiling. The 14 decorations and the  7 equations can be reduced to one. 

\begin{figure}[t]
\begin{center}
\caption{Marked half-hex trapezoid prototile.}
\label{marked half-hex}
\includegraphics[width=1.5in]{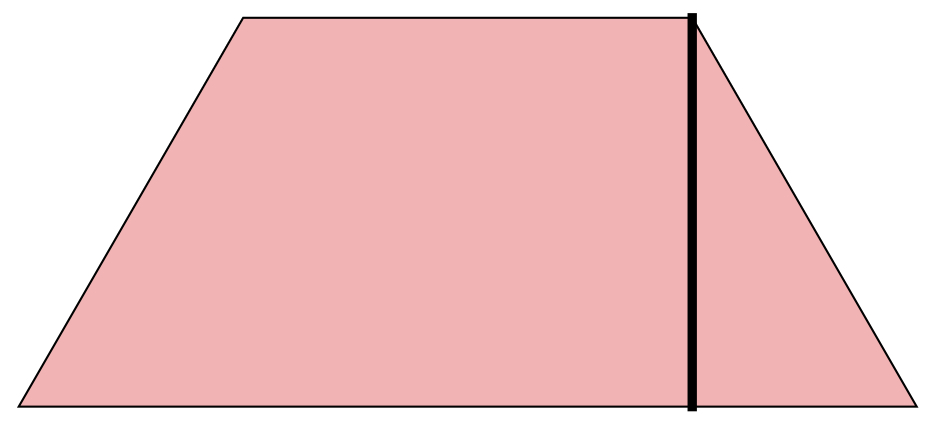}
\end{center}
\end{figure}

Consider the trapezoid shown in fig.\ref{marked half-hex}. It is really only the decoration -- the line segment -- that matters, but the trapezoid helps guide understanding. The trapezoid can be partitioned into four congruent trapezoids similar to the original. This will lead to a self-similar tiling of the plane. The tiling "rule" is that curves described by the repeated line segments must not have endpoints. Thus the curves either start and end at infinity, or they form closed closed curves in compact subsets of the plane.

 \thm The recursion eq(\ref{HLR mono trap recursion}) with $\rho=2$, $\omega=\exp( i \pi /3)$ produces the Taylor-Socolar aperiodic  monotiling: 
\be
\rho \, R = R \cup 
	   \left( -1 + \omega^2 + \omega^4 R   \right) \cup
	 \left( 1 + \omega + \omega^5 \overline{R}   \right) \cup
	  \left( i \sqrt{3}  +  \overline{R}   \right) 
\label{HLR mono trap recursion}
\ee

\begin{figure}[t]
\begin{center}
\caption{A portion of the trapezoid-based monotiling produced by eq(\ref{HLR mono trap recursion}).  The field of trapezoids is shown in the background for reference.  5 iterations.}
\label{HLR mono trap 5 dig}
\includegraphics[width=3in]{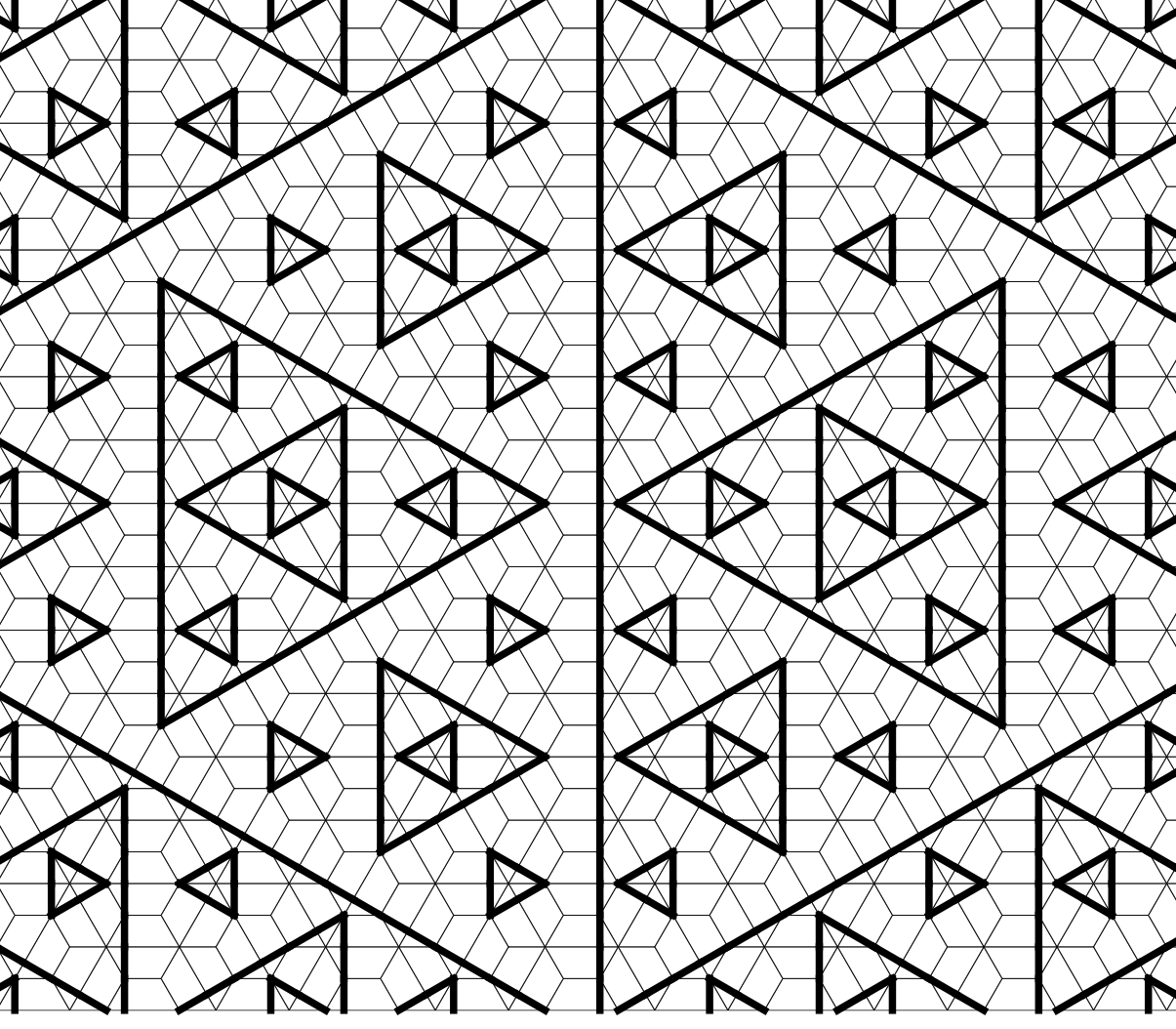}
\end{center}
\end{figure}

\pf Once again,we find the Taylor triangle pattern, a portion of which is shown in fig.\ref{HLR mono trap 5 dig} after 5 iterations. $\Box$\\

Other recursions for the same marked prototile can be written that lead to periodic tilings. It is neither the tile shape, nor the decoration, that uniquely describes aperiodicity. It is the algebraic form of the recursion.

\noindent Boston, MA\\
20141211

\end{document}